# An Efficient Semismooth Newton Based Algorithm for Convex Clustering

Yancheng Yuan,[*] Defeng Sun[†] and Kim-Chuan Toh[‡]

February 20, 2018

**Abstract**

Clustering may be the most fundamental problem in unsupervised learning which is still active in machine learning research because its importance in many applications. Popular methods like K-means, may suffer from instability as they are prone to get stuck in its local minima. Recently, the sum-of-norms (SON) model (also known as clustering path), which is a convex relaxation of hierarchical clustering model, has been proposed in [7] and [5] Although numerical algorithms like ADMM and AMA are proposed to solve convex clustering model [2], it is known to be very challenging to solve large-scale problems. In this paper, we propose a semi-smooth Newton based augmented Lagrangian method for large-scale convex clustering problems. Extensive numerical experiments on both simulated and real data demonstrate that our algorithm is highly efficient and robust for solving large-scale problems. Moreover, the numerical results also show the superior performance and scalability of our algorithm compared to existing first-order methods.

## 1 Introduction

Clustering is one of the most fundamental but important problems in unsupervised learning. Traditional clustering models such as K-means clustering, hierarchical clustering may suffer from instability because of the non-convexity of the model and difficulties in finding a global optimal solution. The clustering results are generally highly dependent on the initialization and the results could be totally different with different initializations. Moreover, these clustering models require the prior knowledge about the number of clusters which are often not available in many real applications. Therefore, in real applications, k-means is typically tried with different cluster numbers and the user will then decide on a suitable value based on his judgment on which computed result agrees best with his domain knowledge or experience. Obviously, such a process could make the clustering results to be subjective.

In order to overcome the above issues, a new clustering algorithm has been proposed recently [7, 5] which is more robust compared to those traditional ones. Let $A \in \mathbb{R}^{d \times n} = [\mathbf{a}_1, \mathbf{a}_2, \cdots, \mathbf{a}_n]$ be a given data matrix with $n$ observations and $d$ features. Convex clustering model for these $n$ observations solves the following convex optimization problem:

$$\min_{X \in \mathbb{R}^{d \times n}} \frac{1}{2} \sum_{i=1}^{n} \|\mathbf{x}_i - \mathbf{a}_i\|^2 + \gamma \sum_{i<j} \|\mathbf{x}_i - \mathbf{x}_j\|_q, \qquad (1)$$

where $\|\cdot\|$ denotes the 2-norm, $q \geq 1$ ensures the convexity and $X = [\mathbf{x}_1, \ldots, \mathbf{x}_n]$. The popular selections of $q$ for convex clustering model (1) are $q \in \{1, 2, \infty\}$. In this paper, we focus on $q = 2$. After solving (1) and obtaining the optimal solution $X^*$, $\mathbf{a}_i$ and $\mathbf{a}_j$ belong to the same cluster if and only if $\mathbf{x}_i^* = \mathbf{x}_j^*$. In other words, $\mathbf{x}_i^*$ is the centroid for observation $\mathbf{a}_i$. The idea behind this model is that if two input observations $\mathbf{a}_i$ and $\mathbf{a}_j$ belong to same cluster, then their corresponding centroids $\mathbf{x}_i$ and $\mathbf{x}_j$ should be the same. The first term in (1) is the fidelity term. The second term is the regularization term to penalize the difference between different centroids and enforce the property that centroids for observations in the same cluster should be identical. From another point of view, we can also regard the second term as $\ell_{1,2}$ norm which can result in group sparsity for the centroids. Therefore, we can expect that there are only a few distinct centroids.

---

[*]Department of Mathematics, National University of Singapore, 10 Lower Kent Ridge Road, Singapore (yuanyancheng@u.nus.edu.sg).

[†]Department of Applied Mathematics, Hong Kong Polytechnic University (defeng.sun@polyu.edu.hk).

[‡]Department of Mathematics, National University of Singapore, 10 Lower Kent Ridge Road, Singapore (mattohkc@nus.edu.sg).

The advantages of convex clustering lie mainly in two aspects. First, since the clustering model (1) is strongly convex, the optimal solution for a given positive $\gamma$ is unique and is more easily obtainable than traditional clustering algorithms like K-means. Second, instead of requiring the prior knowledge of the cluster number, we can generate a clustering path via solving (1) for a sequence of positive values of $\gamma$. To handle cluster recovery for large-scale data sets, various researchers have suggested the following sparse convex clustering model modified from (1)

$$\min_{X \in \mathbb{R}^{d \times n}} \frac{1}{2} \sum_{i=1}^{n} \|\mathbf{x}_i - \mathbf{a}_i\|^2 + \gamma \sum_{(i,j) \in \mathcal{E}} w_{ij} \|\mathbf{x}_i - \mathbf{x}_j\|_q, \qquad (2)$$

where the edge set $\mathcal{E} = \{(i, j) \mid \text{if } j \text{ is among } i\text{'s } k\text{-nearest neighbors}\}$ and $w_{ij} = \exp(-\phi \|\mathbf{a}_i - \mathbf{a}_j\|^2)$ for $(i, j) \in \mathcal{E}$. We can regard the original convex clustering model (1) as a special case if we take $k = n$ and $\phi = 0$.

The advantages just mentioned and the success of the convex model (1) in recovering clusters in many examples with well selected value of $\gamma$ have motivated researchers to provide theoretical guarantees on the cluster recovery property of (1). However, the first theoretical result on cluster recovery established in [14] is valid for only two clusters. It showed that the model (1) can recover the two clusters perfectly if the data points are drawn from two cubes and the distance between these two cubes are large enough. In 2015, Tan et al. [13] analyzed the statistical properties of (1). Recently, Panahi et al. [9] provided theoretical results for general $k$ clusters case under relatively mild sufficient conditions, thus laying the theoretical foundation for convex clustering. However, the conditions provided in the theoretical analysis are usually not checkable before we find the right clusters and the perfect value of $\gamma$ is unknown. In practice, we need to try a sequence of values of $\gamma$ to generate a clustering path.

The challenges for the convex model (1) to obtain meaningful cluster recovery is then to solve it efficiently for a range of values of $\gamma$. Lindsten [7] used an off-the-shelf convex solver, CVX to generate the solution path. However, Hocking et al. [5] realized that CVX is competitive only when the problem size is small and it does not scale well when increases the number of data points. Thus the paper introduced three algorithms for different regularizers. Recently, some new algorithms have been proposed to solve this problem. Chi & Lange [2] introduced ADMM and AMA to solve (1). However, based on our experiments, both algorithms may still encounter scalability issues, albeit less severe than CVX. Futhermore, the efficiency of these two algorithms are sensitive to the parameter value $\gamma$. This is not favorable since we need to solve (1) with $\gamma$ in a relative large range to generate the clustering path. In [9], Panahi et al. proposed a stochastic splitting algorithm to solve (1) in an attempt to resolve the aforementioned scalability issues. Although this stochastic approach scales well with the problem scale ($n$ in (1)), the convergence rate shown in [9] is rather weak in that it requires at least $l \geq n^4/\varepsilon$ iterations to generate a solution $X^l$ such that $\|X^l - X^*\|^2 \leq \varepsilon$ is satisfied with high probability. (Here and below, $\|\cdot\|$ is also used to denote the Frobenius norm of a matrix.) Moreover, because the error estimate is given in the sense of high probability, it is difficult to design appropriate stopping condition for the algorithm in practice.

As the reader may observe, all the existing algorithms are purely first-order methods that do not use any second-order information underlying the convex clustering model. In this paper, we design and analyse a deterministic second-order algorithm, called a semi-smooth Newton augmented Lagrangian method, to solve the convex clustering model. The algorithm is not only proven to be theoretically efficient but it is also demonstrated to be practically highly efficient and robust.

## 2 Related Work

The convex clustering model (1) was proposed by [7] and [5] in 2011. Lindsten et al. [7] used the off-the shelf solver, CVX, to solve it, however, Hocking [5] proposed three different algorithms to solve (1) with $q \in \{1, 2, \infty\}$, respectively. In particular, Hocking used the subgradient method for solving (1) with $\ell_2$ norm regularization. It's until 2013, Chi et al. [2] introduced two unified frameworks, i.e. ADMM and AMA, to solve (1) and (2) with $q \in \{1, 2, \infty\}$. In 2016, [6] proposed a highly efficient semi-smooth Newton ALM framework to solve Lasso and fused Lasso problem, which inspires us to propose the algorithmic frame for convex clustering model in this paper.

## 3 A semi-smooth Newton-CG augmented Lagrangian method with fast linear convergence

In this section, we introduce a fast convergent augmented Lagrangian method for solving the convex clustering model (2). Before that, we introduce some preliminaries and notation.

### 3.1 Preliminaries and Notation

For a given undirected graph $G = (\mathcal{V}, \mathcal{E})$ with $n$ vertices and edges defined in $\mathcal{E}$, we define the symmetric adjacency matrix $J \in \mathbb{R}^{n \times n}$ with entries

$$J_{ji} = J_{ij} = \left\{ \begin{array}{ll} 1 & if \ (i,j) \in \mathcal{E}, \\ 0 & otherwise. \end{array} \right.$$

We also define the node-arc incidence matrices $\mathcal{J}, \bar{\mathcal{J}} \in \mathbb{R}^{n \times |\mathcal{E}|}$ as

$$\mathcal{J}_k^{(i,j)} = \left\{ \begin{array}{ll} 1 & if \ k = i, \\ 0 & otherwise. \end{array} \right. \tag{3}$$

$$\bar{\mathcal{J}}_k^{(i,j)} = \begin{cases} 1 & if\ k = j, \\ 0 & otherwise. \end{cases} \quad (4)$$

Where $\mathcal{J}_k^{(i,j)}$ and $\bar{\mathcal{J}}_k^{(i,j)}$ is the k-th entry of the $(i,j)-th$ column.

**Proposition 1.** *With matrices $J, \mathcal{J}, \bar{\mathcal{J}}$ defined above, we have the following results*

$$\mathcal{J}\mathcal{J}^T + \bar{\mathcal{J}}\bar{\mathcal{J}}^T = diag(J\mathbf{e}), \quad \mathcal{J}\bar{\mathcal{J}}^T + \bar{\mathcal{J}}\mathcal{J}^T = J.$$

*Where $\mathbf{e} \in \mathbb{R}^n$ is the column vector of all ones.*

Now, for given variables $X \in \mathbb{R}^{d \times n}$, $Z \in \mathbb{R}^{d \times |\mathcal{E}|}$ and undirected graph $G = (\mathcal{V}, \mathcal{E})$, we define a linear operator $\mathcal{B}: \mathbb{R}^{d \times n} \to \mathbb{R}^{d \times |\mathcal{E}|}$ and its adjoint operator $\mathcal{B}^*(Z): \mathbb{R}^{d \times |\mathcal{E}|} \to \mathbb{R}^{d \times n}$, respectively, by

$$\mathcal{B}(X) = [(\mathbf{x}_i - \mathbf{x}_j)]_{(i,j) \in \mathcal{E}} = X(\mathcal{J} - \bar{\mathcal{J}}),$$

and

$$\mathcal{B}^*(Z) = Z(\mathcal{J}^T - \bar{\mathcal{J}}^T).$$

Thus, by Proposition (1), we have

$$\mathcal{B}^*(\mathcal{B}(X)) = X(\mathcal{J}\mathcal{J}^T - \mathcal{J}\bar{\mathcal{J}}^T - \bar{\mathcal{J}}\mathcal{J}^T + \bar{\mathcal{J}}\bar{\mathcal{J}}^T) = XL_J, \quad (5)$$

where $L_J = diag(J\mathbf{e}) - J \in \mathbb{R}^{n \times n}$ is the unweighted Laplacian matrix associated with the undirected graph $G$.

For a convex function $p: \mathcal{X} \to \mathcal{Y}$, which is proper and closed, the proximal mapping $\text{Prox}_{tp}(x)$ for $p$ at any $x \in \mathcal{X}$ with $t > 0$ is defined by

$$\text{Prox}_{tp}(x) = \arg\min_{u \in \mathcal{X}} \{tp(u) + \frac{1}{2}\|u - x\|^2\}. \quad (6)$$

It is well known that proximal mappings are important for designing optimization algorithms and they have been well studied. The proximal mappings for many common used functions in fact have closed form formulas. Here, we summarize those that are related to this paper. Note that $\mathcal{P}_C$

Table 1: Proximal maps for selected functions

| $p(\cdot)$ | $\text{Prox}_{tp}(\mathbf{x})$ | Comment |
|---|---|---|
| $\|\cdot\|_1$ | $\left[1 - \frac{t}{|\mathbf{x}_l|}\right]_+ \mathbf{x}_l$ | Elementwise soft-thresholding |
| $\|\cdot\|_2$ | $\left[1 - \frac{t}{\|\mathbf{x}\|_2}\right]_+ \mathbf{x}$ | Blockwise soft-thresholding |
| $\|\cdot\|_\infty$ | $\mathbf{x} - \mathcal{P}_{t\mathcal{S}}(\mathbf{x})$ | $\mathcal{S}$ is the unit $\ell_1$-ball |

denotes the projection onto a given closed convex set $C$.

In this paper, we will often make use of the following Moreau identity

$$\text{Prox}_{tp}(x) + t\text{Prox}_{p^*/t}(x/t) = x,$$

where $t > 0$ and $p^*$ is the conjugate function of $p$.

## 3.2 Duality and Optimality Conditions

In this section, we will derive the dual problem of (2) and the KKT conditions. First, we write (2) equivalently as the following compact form

$$(P) \quad \min_{X,U} \left\{ \frac{1}{2}\|X - A\|^2 + p(U) \mid \mathcal{B}(X) - U = 0 \right\},$$

where $p(U) = \gamma \sum_{(i,j) \in \mathcal{E}} w_{ij} \|U^{(i,j)}\|$. Here $U^{(i,j)}$ denotes the $(i,j)$-th column of $U \in \mathbb{R}^{d \times |\mathcal{E}|}$. The dual problem is given by

$$(D) \quad \max_{V,Z} \left\{ \langle A, V \rangle - \frac{1}{2}\|V\|^2 \mid \mathcal{B}^*(Z) - V = 0, Z \in \Omega \right\},$$

where $\Omega = \{Z \in \mathbb{R}^{d \times |\mathcal{E}|} \mid \|Z^{(i,j)}\| \leq w_{ij}, \forall (i,j) \in \mathcal{E}\}$.

Now, denote by $l$ the Lagrangian function for $(P)$:

$$l(X, U; Z) = \frac{1}{2}\|X - A\|^2 + p(U) + \langle Z, \mathcal{B}(X) - U \rangle. \tag{7}$$

Furthermore, given $\sigma > 0$, the augmented Lagrangian function associated with $(P)$ is given by

$$\mathcal{L}_\sigma(X, U; Z) = l(X, U; Z) + \frac{\sigma}{2}\|\mathcal{B}(X) - U\|^2.$$

The KKT conditions for the primal-dual problems are given by

$$(KKT) \quad \begin{cases} V + X - A &= \mathbf{0}, \\ U - \text{Prox}_p(U + Z) &= \mathbf{0}, \\ \mathcal{B}(X) - U &= \mathbf{0}, \\ \mathcal{B}^*(Z) - V &= \mathbf{0}. \end{cases}$$

## 3.3 A Semi-Smooth Newton-CG augmented Lagrangian Method for Solving (2)

In this section, we will design an inexact augmented Lagrangian method for solving the primal problem $(P)$ but it will also solve $(D)$ as a byproduct. Since a semismooth Newton-CG method will be used to solve the subproblems involved in the method, we call our algorithm a semismooth Newton-CG augmented Lagrangian method (SSNAL in short). The algorithm for solving $(P)$ is shown in **Algorithm** 1. To ensure the convergence of the inexact augmented Lagrangian method in **Algorithm** 1, we need the following stopping criterion for solving the subproblem (8) in each iteration:

$$(A) \quad \text{dist}(0, \partial \Phi_k(X^{k+1}, U^{k+1})) \leq \epsilon_k / \max\{1, \sqrt{\sigma_k}\},$$

where $\{\epsilon_k\}$ is a given summable sequence of nonnegative numbers.

## 3.4 Solving the subproblem (8)

The inexact augmented Lagrangian method (ALM) is a well studied and efficient theoretical algorithmic framework for composite convex optimization problems. The key challenge in making the ALM efficient numerically is solving the associated subproblem (8) in each iteration efficiently to a

**Algorithm 1** S<small>SNAL</small> for $(P)$

**Initialization**: Choose $(X^0, U^0, Z^0) \in \mathbb{R}^{d\times n} \times \mathbb{R}^{d\times|\mathcal{E}|} \times \mathbb{R}^{d\times|\mathcal{E}|}$, $\sigma_0 > 0$ and a summable nonnegative error tolerance sequence $\{\epsilon_k\}$.

**repeat**

  **Step 1**. Compute
$$(X^{k+1}, U^{k+1}) \approx \arg\min\{\Phi_k(X,U) = \mathcal{L}_{\sigma_k}(X,U;Z^k)\} \tag{8}$$
  to satisfy the stopping condition (A).

  **Step 2**. Compute
$$Z^{k+1} = Z^k + \sigma_k(\mathcal{B}(X^{k+1}) - U^{k+1}).$$

  **Step 3**. Update $\sigma_{k+1} \uparrow \sigma_\infty \leq \infty$.

**until** Stopping criterion is satisfied.

---

required accuracy. In this paper, we will design a semismooth Newton-CG method to solve (8). We will establish its quadratic convergence and derive sophisticated numerical techniques to solve the semismooth Newton equations very efficiently by exploiting the underlying second-order structured sparsity in the subproblem.

For given $\sigma$ and $\tilde{Z}$, the subproblem (8) in each iteration has the following form:
$$\min_{X\in\mathbb{R}^{d\times n}, U\in\mathbb{R}^{d\times|\mathcal{E}|}} \Phi(X,U) := \mathcal{L}_\sigma(X,U;\tilde{Z}). \tag{9}$$

Since $\Phi(X,U)$ is a strongly convex function, the level set $\mathcal{L}_\alpha\{(X,U)|\Phi(X,U) \leq \alpha\}$ is a closed and bounded convex set for any $\alpha \in \mathbb{R}$. Moreover, (9) admits a unique optimal solution denoted as $(\bar{X},\bar{U})$. Now, for any $X$, denote
$$\begin{aligned}\phi(X) &:= \inf_U \Phi(X,U) \\ &= \tfrac{1}{2}\|X-A\|^2 + p(\text{Prox}_{p/\sigma}(\mathcal{B}(X) + \sigma^{-1}\tilde{Z})) \\ &\quad + \tfrac{1}{2\sigma}\|\text{Prox}_{\sigma p^*}(\sigma\mathcal{B}(X) + \tilde{Z})\|^2 - \tfrac{1}{2\sigma}\|\tilde{Z}\|^2.\end{aligned}$$

Therefore, we can compute $(\bar{X},\bar{U}) = \arg\min \Phi(X,U)$ via
$$\bar{X} = \arg\min \phi(X), \quad \bar{U} = \text{Prox}_{p/\sigma}(\mathcal{B}(\bar{X}) + \sigma^{-1}\tilde{Z}).$$

Since $\phi(\cdot)$ is strongly convex and continuously differentiable on $\mathbb{R}^{d\times n}$ with
$$\nabla\phi(X) = X - A + \mathcal{B}^*(\text{Prox}_{\sigma p^*}(\sigma\mathcal{B}(X) + \tilde{Z})),$$
we know that $\bar{X}$ can be obtained by solving the following nonsmooth equation
$$\nabla\phi(X) = 0. \tag{10}$$

It is well known that Newton's method is the best method to solve nonlinear equations if it can be implemented efficiently. However, Newton's method requires the smoothness of $\nabla\phi(X)$ which is not the case in our problem. This motivates us to develop a semismooth Newton method to solve the nonsmooth equation (10). Before we present our semismooth Newton method, we introduce the following definition of semismoothness.

**Definition 1.** *(Semismoothness). Let $F : \mathcal{O} \subseteq \mathcal{X} \to \mathcal{Y}$ be a locally Lipschitz continuous function on the open set $\mathcal{O}$. $F$ is said to be semismooth at $x \in \mathcal{O}$ if $F$ is directionally differentiable at $x$ and for any $V \in \partial F(x + \Delta x)$ with $\Delta x \to 0$,*

$$F(x + \Delta x) - F(x) - V\Delta x = o(\|\Delta x\|).$$

*$F$ is said to be strongly semismooth at $x$ if $F$ is semismooth at $x$ and*

$$F(x + \Delta x) - F(x) - V\Delta x = O(\|\Delta x\|^2).$$

*$F$ is said to be a semismooth (respectively, strongly semismooth) function on $\mathcal{O}$ if it is semismooth (respectively, strongly semismooth) everywhere in $\mathcal{O}$.*

**Lemma 1.** *For any $t > 0$, the proximal mapping $\mathrm{Prox}_{t\|\cdot\|_2}$ is strongly semismooth.*

Now, we derive the generalized Jacobian of the locally lipschitz continuous function $\nabla \phi(\cdot)$. For any given $X \in \mathbb{R}^{d \times n}$, the following set-valued map is well defined:

$$\begin{aligned}
\hat{\partial}^2 \phi(X) &:= \{\mathcal{I} + \sigma \mathcal{B}^* \mathcal{V} \mathcal{B} | \mathcal{V} \in \partial \mathrm{Prox}_{\sigma p^*}(\tilde{Z} + \sigma \mathcal{B}(X))\} \\
&= \{\mathcal{I} + \sigma \mathcal{B}^* (\mathcal{I} - \mathcal{P})\mathcal{B} | \mathcal{P} \in \partial \mathrm{Prox}_{\sigma^{-1}p}(\sigma^{-1}\tilde{Z} + \mathcal{B}(X))\},
\end{aligned}$$

where $\partial \mathrm{Prox}_{\sigma p^*}(\tilde{Z} + \sigma \mathcal{B}(X))$ and $\partial \mathrm{Prox}_{p/\sigma}(\sigma^{-1}\tilde{Z} + \mathcal{B}(X))$ are the Clarke subdifferential of the Lipschitz continuous mapping $\mathrm{Prox}_{\sigma p^*}$ and $\mathrm{Prox}_{p/\sigma}(\cdot)$ at $\tilde{Z} + \sigma \mathcal{B}(X)$ and $\sigma^{-1}\tilde{Z} + \mathcal{B}(X)$, respectively. Note that from [3] [p.75] and [4] [Example 2.5], we have

$$\partial^2 \phi(X)(d) = \hat{\partial}^2 \phi(X)(d), \quad \forall d \in \mathbb{R}^{d \times n},$$

where $\partial^2 \phi(X)$ is the generalized Hessian of $\phi$ at $X$. Thus we may use $\hat{\partial}^2 \phi(X)$ as the surrogate for $\partial^2 \phi(X)$. Since $\mathcal{I} - \mathcal{P} = \mathcal{V} \in \partial \mathrm{Prox}_{\sigma p^*}(\cdot)$ is symmetric and positive semdefinite. The elements in $\hat{\partial}^2 \phi(X)$ are positive definite, which will guarantee that (11) is well defined.

Now, we can present our semismooth Newton-CG (SSNCG) method for solving (10) and we could expect to get a fast superlinear or even quadratic convergence.

---

**Algorithm 2** SSNCG for (10)

---

**Initialization**: Given $X^0 \in \mathbb{R}^{d \times n}$, $\mu \in (0, 1/2)$, $\tau \in (0, 1]$, and $\bar{\eta}, \delta \in (0, 1)$. For $j = 0, 1, \ldots$
**repeat**
   **Step 1**. Apply the conjugate gradient (CG) method to find an approximate solution $d^j \in \mathbb{R}^{d \times n}$ to

$$\mathcal{V}_j(d) \approx -\nabla \phi(X^j) \qquad (11)$$

with $\mathcal{V}_j \in \hat{\partial}^2 \phi(X^j)$ defined in (10), such that $\|\mathcal{V}_j(d^j) + \nabla \phi(X^j)\| \leq \min(\bar{\eta}, \|\nabla \phi(X^j)\|^{1+\tau})$.
   **Step 2**. (Line Search) Set $\alpha_j = \delta^{m_j}$, where $m_j$ is the first nonnegative integer $m$ for which

$$\phi(X^j + \delta^m d^j) \leq \phi(X^j) + \mu \delta^m \langle \nabla \phi(X^j), d^j \rangle.$$

   **Step 3**. Set $X^{j+1} = X^j + \alpha_j d^j$.
**until** Stopping criterion is satisfied.

---

## 3.5 Using the conjugate gradient method to solve (11)

In this section, we will discuss how to solve the very large symmetric positive definite linear system (11) to compute the Newton direction efficiently. As the coefficient matrix $\mathcal{V}_j$ in (11) is expensive to compute and factorize, we will adopt the conjugate gradient (CG) method to solve it. The computational cost for CG is highly dependent on the cost for computing the matrix-vector product $\mathcal{V}_j(\tilde{d})$ for any given $\tilde{d} \in \mathbb{R}^{d \times n}$. Thus we will need to analyse how this product can be computed efficiently. We start by choosing $\mathcal{P} \in \partial \text{Prox}_{p/\sigma}(\mathcal{B}(X^j) + \sigma^{-1}\tilde{Z})$ explicitly. Let $D := \mathcal{B}(X^j) + \sigma^{-1}\tilde{Z}$. We can take $\mathcal{P} : \mathbb{R}^{d \times |\mathcal{E}|} \to \mathbb{R}^{d \times |\mathcal{E}|}$ that is defined by

$$(\mathcal{P}(\tilde{U}))^{(i,j)} = \begin{cases} 0 & if\ \alpha_{ij} \geq 1, \\ \alpha_{ij} \frac{\langle D^{(i,j)}, \tilde{U}^{(i,j)} \rangle}{\|D^{(i,j)}\|^2} D^{(i,j)} \\ + (1 - \alpha_{ij})\tilde{U}^{(i,j)} & if\ \alpha_{ij} < 1, \end{cases}$$

where for $(i, j) \in \mathcal{E}$,

$$\alpha_{ij} = \begin{cases} \frac{\sigma^{-1} w_{ij}}{\|D^{(i,j)}\|} & \text{if } \|D^{(i,j)}\| > 0, \\ \infty & \text{if } \|D^{(i,j)}\| = 0. \end{cases}$$

Note that for given $D \in \mathbb{R}^{d \times |\mathcal{E}|}$, the cost for computing $\alpha$ is $2(d+1)|\mathcal{E}|$ arithmetic operations.

For convenience, denote

$$\widehat{\mathcal{E}} = \{(i,j) \in \mathcal{E} \mid \alpha_{ij} < 1\}.$$

To compute $\mathcal{V}_j(X) = X + \sigma \mathcal{B}^* \mathcal{B}(X) - \mathcal{B}^* \mathcal{P} \mathcal{B}(X)$ efficiently for a given $X \in \mathbb{R}^{d \times n}$, we need the efficient computation of $\mathcal{B}^* \mathcal{P} \mathcal{B}(X)$ as presented in the following proposition.

**Proposition 2.** *Consider the symmetric matrix $M \in \mathbb{R}^{n \times n}$ defined by $M_{ij} = \alpha_{ij}$ if $(i,j) \in \widehat{\mathcal{E}}$ and $M_{ij} = 0$ otherwise. Let $Y = [M_{ij}(\mathbf{x}_i - \mathbf{x}_j)]_{(i,j)\in\mathcal{E}} = X(\mathcal{M} - \bar{\mathcal{M}})$, where $\mathcal{M}$ and $\bar{\mathcal{M}}$ are defined similarly as in (3) and (4), respectively, for the matrix $M$. Then we have*

$$\mathcal{B}^*(Y) = X L_M,$$

*where $L_M$ is the Laplacian matrix associated with $M$. The cost of computing the result is $O(d|\widehat{\mathcal{E}}|)$ arithmetic operations.*

Now, define $\rho \in \mathbb{R}^{|\mathcal{E}|}$ by

$$\rho_{(i,j)} := \begin{cases} \frac{1}{\|D^{(i,j)}\|^2} \langle D^{(i,j)}, \mathbf{x}_i - \mathbf{x}_j \rangle, & \text{if } (i,j) \in \widehat{\mathcal{E}}, \\ 0, & \text{otherwise.} \end{cases}$$

For given $D \in \mathbb{R}^{d \times |\mathcal{E}|}$, the cost for computing $\rho$ is $O(d|\widehat{\mathcal{E}}|)$ arithmetic operations. Let $W_{(i,j)} = \frac{\langle D^{(i,j)}, \mathbf{x}_i - \mathbf{x}_j \rangle}{\|D^{(ij)}\|^2} D^{(i,j)} = \rho_{(i,j)} D^{(i,j)}$. Then,

$$\mathcal{B}^*(W) = W(\mathcal{J}^T - \bar{\mathcal{J}}^T) = D \text{diag}(\rho)(\mathcal{J}^T - \bar{\mathcal{J}}^T).$$

Thus, the computational cost for $\mathcal{B}^* \mathcal{P} \mathcal{B}(X)$ in total is $O(d|\widehat{\mathcal{E}}|)$.

Observe that $\alpha_{ij} < 1$ means $j$ is in $i$'s nearest $k$ neighbors but does not belong to the same cluster as $i$, so $|\widehat{\mathcal{E}}|$ should be much smaller than $|\mathcal{E}|$. On the other hand, for $\alpha_{ij} \geq 1$, it means that points $i$ and $j$ are in the same cluster. From this analysis, we can expect most of the columns of the matrix $\mathcal{P}(\mathcal{B}(X))$ to be zero. We call such a property inherited from the generalized Hessian of $\phi(\cdot)$ at $X$ as the **second-order sparsity.** It is because of this important property that we can able to compute $\mathcal{B}^* \mathcal{P} \mathcal{B}(X)$ at a very low cost.

## 3.6 Convergence Results

In this section, we will establish the convergence rate result for both SSNAL and SSNCG under mild assumptions. First, we present the following global convergence result of our proposed Algorithm SSNAL.

**Theorem 1.** *Let $\{(X^k, U^k, Z^k)\}$ be the sequence generated by Algorithm 1 with stopping criterion (A). Then the sequence $\{X^k\}$ is bounded and converges to the unique optimal solution of $(P)$, and $\|\mathcal{B}(X^k) - U^k\|$ converges to $0$. In addition, $\{Z^k\}$ is also bounded and converges to the an optimal solution $Z^* \in \Omega$ of $(D)$.*

The above convergence theorem can be obtained from [11, 12] without much difficulties. Next, we state the convergence property for the semismooth Newton algorithm SSNCG used to solve the subproblems in Algorithm 1.

**Theorem 2.** *Assume that $Prox_{tp}(\cdot)$ is strongly semismooth on $int(dom(p))$. Let the sequence $\{X^j\}$ be generated by Algorithm SSNCG. Then $\{X^j\}$ converges to the unique solution $\bar{X}$ of the problem in (10) and*
$$\|X^{j+1} - \bar{X}\| = O(\|X^j - \bar{X}\|^{1+\tau}),$$
*where $\tau \in (0, 1]$ is a given constant in the algorithm.*

The proof of this theorem could be found in the supplementary material. Note that by Lemma 1, the strongly semismoothness assumption holds true for our model (2).

## 3.7 Generating Initial Points

In our implementation, we use the following inexact alternating direction method of multipliers (IADMM) developed in [1] to generate an initial point to warm-start SSNAL[1]. Note that in Step 1, $X^{k+1}$ is a computed solution for the following large linear system of equations:
$$(I_n + \sigma \mathcal{B}^* \mathcal{B})X = R^k \underset{(5)}{\iff} (I_n + \sigma L_J)X^T = (R^k)^T.$$

To compute $X^{k+1}$, we apply the conjugate gradient method to solve the above linear system.

## 4 Numerical Experiments

In this section, we show the superior performance of our proposed algorithm SSNAL on both simulated and real datasets, comparing to the popular algorithms such as ADMM and AMA which are proposed in [2]. In particular, we will focus on the efficiency, scalability, and robustness of our algorithm for different values of $\gamma$. Also, we will show the performance of our algorithm on large datasets and unbalanced data. Previous research on scalability and performance of (2) on unbalanced datasets is limited. The problem sizes of the instances tested in [2] and other related papers are only several hundreds ($n \leq 500$ in [2], $n \leq 600$ in [9]), which are not large enough to clearly demonstrate the scalability of the algorithms. In this paper, we will show numerical results

---
[1] With the global convergence result stated in Theorem 1, the performance of SSNAL does not sensitively depend on the initial points, but it is still helpful if we can choose a good one.

---

**Algorithm 3** IADMM for $(P)$

---

**Initialization**: Choose $\sigma > 0$, $(X^0, U^0, Z^0) \in \mathbb{R}^{d \times n} \times \mathbb{R}^{d \times |\mathcal{E}|} \times \mathbb{R}^{d \times |\mathcal{E}|}$, and a summable nonnegative error tolerance sequence $\{\epsilon_k\}$. For $k = 0, 1, \ldots$,

**repeat**

  **Step 1**. Let $R^k = A + \sigma \mathcal{B}^*(U^k - \sigma^{-1} Z^k)$. Compute

$$\begin{aligned} X^{k+1} &\approx \arg\min_X \{\mathcal{L}_\sigma(X, U^k; Z^k)\}, \\ U^{k+1} &= \arg\min_U \{\mathcal{L}_\sigma(X^{k+1}, U; Z^k)\}. \end{aligned}$$

  where $X^{k+1}$ is an inexact solution satisfying the accuracy requirement that $\|(I_n + \sigma \mathcal{B}^* \mathcal{B}) X^{k+1} - R^k\| \leq \epsilon_k$.

  **Step 2**. Compute

$$Z^{k+1} = Z^k + \tau \sigma_k (\mathcal{B}(X^{k+1}) - U^{k+1}),$$

  where $\tau \in (0, \frac{1+\sqrt{5}}{2})$ is typically chosen to be 1.618.

**until** the stopping criterion is satisfied.

---

for $n$ up to 20000. Also, we will analyze the sensitivity of the computational efficiency of SSNAL and AMA, with respect to different choices of hyper-parameters in (2), such as $k$ (number of nearest neighbors) and $\gamma$.

Our attention in this paper focuses on solving (2) with $q = 2$ since the rational invariance of the $\ell_2$ norm makes it a robust choice in practice. Also, this case is more challenging than $q = 1$ or $q = \infty$.[2] As the results reported in [2] have been regarded as the benchmark for the convex clustering model (2), we will compare our algorithm with the open source software CVXCLUSTR[3] in [2], which is an R package with key functions written in C. We write our code in MATLAB without any dedicated C functions. All our computational results are obtained from a desktop having 16 cores with 32 Intel Xeon E5-2650 processors at 2.6 GHz and 64 GB memory.

In our implementation, we stop our algorithm by the following relative KKT residual:

$$\max\{\eta_P, \eta_D, \eta\} \leq \epsilon$$

where

$$\eta_P = \frac{\|\mathcal{B}X - U\|}{1 + \|U\|},$$

$$\eta_D = \frac{\sum_{(i,j) \in \mathcal{E}} \max\{0, \|Z^{(i,j)}\|_2 - w_{ij}\}}{1 + \|A\|},$$

$$\eta = \frac{\|\mathcal{B}^*(Z) + X - A\| + \|U - \text{Prox}_p(U + Z)\|}{1 + \|A\| + \|U\|},$$

and $\epsilon > 0$ is a given tolerance. Since the numerical results reported in [2] have demonstrated the superior performance of AMA over ADMM, we will mainly compare our proposed algorithm with AMA. We note that CVXCLUSTR does not use the relative KKT residual as its stopping criterion but used the duality gap in AMA and $\max\{\eta_P, \eta_D\} \leq \epsilon$ in ADMM. To make a fair comparison, we first solve (2) using SSNAL with a given tolerance $\epsilon$, and denote the primal objective value obtained as $P_{Ssnal}$. Then, we run AMA in CVXCLUSTR and stop it as soon as the computed primal objective function value ($P_{AMA}$) is close enough to $P_{Ssnal}$, i.e.,

$$P_{AMA} - P_{Ssnal} \leq 10^{-6} P_{Ssnal}. \tag{12}$$

---

[2] Our algorithm can be generalized to solve (2) with $q = 1$ and $q = \infty$ easily.
[3] https://cran.r-project.org/web/packages/cvxclustr/index.html

We note that since (2) is an unconstrained problem, the quality of the computed solutions can directly be compared based on the objective function values. We also stop AMA if the maximum of $10^5$ iterations is reached.

In our experiments, we will set $\epsilon = 10^{-6}$ unless specified otherwise. When we generate the clustering path for the first parameter value of $\gamma$, we first run the IADMM introduced in Algorithm 3 for 100 iterations to generate an initial point, then we use SSNAL to solve (2). After that, we use the previously computed optimal solution for the lastest $\gamma$ as the initial point to warm-start SSNAL for solving the problem corresponding to the next $\gamma$. The same strategy is used in CVXCLUSTR.

## 4.1 Simulated Data

In this section, we show the performance of our algorithm SSNAL on two simulated datasets: Two Half Moon and Unbalanced Gaussian [10]. We compare our SSNAL with the AMA in [2] on different problem scales. The numerical results in Table 2 show the superior performance of SSNAL. We also visualize some selected recovery results for Two Half moon and Unbalanced Gaussian in Figure 1.

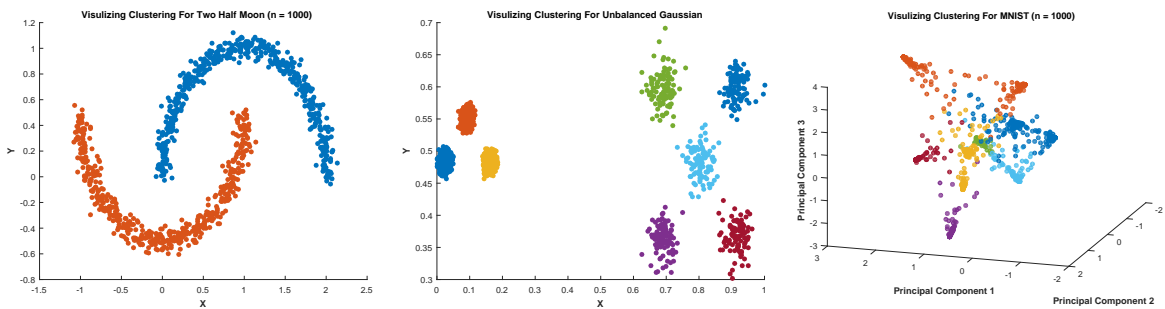

Figure 1: Selected recovery results by model (2) with $\ell_2$ norm. Left: recovery result for Two Half Moon Dataset with $n = 1000$, $k = 20$, $\gamma = 5$. Middle: recovery result for Unbalanced Gaussian Dataset with $n = 6500$, $k = 10$, $\gamma = 1$. Right: recovery result for subset of MNIST with $n = 1000$, $\gamma = 1$.

**Two Half Moon**

The simulated data of two interlocking half moons in $\mathbb{R}^2$ is one of the most popular simulated datasets in clustering. Here we compare the computational time between our proposed SSNAL and AMA on this dataset with different problem scales. We note that AMA could not satisfy the stopping criteria (12) within 100000 iterations when $n$ is large. In the experiments, we choose $k = 10$, $\phi = 0.5$ (for the weights $w_{ij}$) and $\gamma \in [0.2 : 0.2 : 10]$ to generate the clustering path. After generating the clustering path with Ssnal, we repeat the experiments using the same pre-stored primal objective values and stop the AMA using the criterion (12). We report the average time for solving each problem (50 in total) in Table 4.1. Observe that our Ssnal can be more than 50 times faster than AMA.

Table 2: Computation Time (in seconds). Comparison on Two Half Moon. (— means maximum iterations is reached)

| $n$ | 200 | 500 | 1000 | 2000 | 5000 | 10000 |
|---|---|---|---|---|---|---|
| AMA | 0.41 | 4.43 | 28.27 | 78.36 | — | — |
| Ssnal | 0.11 | 0.15 | 0.51 | 1.63 | 7.69 | 20.96 |

**Unbalanced Gaussian**

Next, we show the performance of Ssnal and AMA on the Unbalanced Gaussian data [10]. We can see from Figure 1 that the convex clustering model (2) can recover the cluster assignments perfectly with well chosen parameters. In this experiment, we solve (2) with $k = 10$, $\phi = 0.5$ and $\gamma \in [0.2 : 0.2 : 2]$. For this dataset, we have scaled it so that each entry is in the interval $[0, 1]$.

In experiment, we find that AMA have great difficulty in reaching the stopping criterion (12). We summarize some selected results in Table 3, wherein we report the computation time and iterations for both AMA and Ssnal. Since the main cost for Ssnal is in the inner iterations of Ssncg, so we show the iterations of Ssncg for comparison. The computed primal objective values and other results are reported in the supplementary material.

Table 3: Numerical Results on Unbalanced Gaussian Dataset. $k = 10$, $\phi = 0.5$.

| $\gamma$ | 0.2 | 0.4 | 0.6 | 0.8 | 1.0 |
|---|---|---|---|---|---|
| $t_{\text{AMA}}$ | 264.54 | 256.21 | 260.06 | 262.16 | 263.27 |
| $t_{\text{Ssnal}}$ | **1.15** | **0.57** | **0.67** | **0.66** | **0.86** |
| Iter$_{\text{AMA}}$ | 100000 | 97560 | 97333 | 100000 | 100000 |
| Iter$_{\text{Ssncg}}$ | **23** | **21** | **24** | **24** | **27** |

## 4.2 Real Data

In this section, we compare the performance of our proposed Ssnal with AMA on some real datasets, namely, MNIST[4], Fisher Iris[5], WINE[6], Yale Face B(10Train subset)[7]. For real datasets, a preprocessing step is sometimes necessary to transform the data to one whose features are meaningful for clustering. Thus, for a subset of MNIST (we selected a subset because AMA cannot handle the whole dataset), we first apply the preprocessing method described in [8]. Then we apply the model (2) on the preprocessed data. We summarize the comparison results between SSNAL and AMA on real datasets in Table 4.

---

[4] http://yann.lecun.com/exdb/mnist/

[5] https://archive.ics.uci.edu/ml/datasets/iris

[6] https://archive.ics.uci.edu/ml/datasets/wine

[7] http://www.cad.zju.edu.cn/home/dengcai/Data/FaceData.html

Table 4: Computation Time Comparison on Real Data.

| Dataset | $d$ | $n$ | AMA(s) | SSNAL(s) |
|---|---|---|---|---|
| MNIST | 10 | 1000 | 79.48 | **1.54** |
| Fisher Iris | 4 | 150 | 0.58 | **0.16** |
| WINE | 13 | 178 | 2.62 | **0.19** |
| Yale Face B | 1024 | 760 | 211.36 | **52.71** |

### 4.3 Sensitivity with Different $\gamma$

In order to generate a clustering path for a given dataset, we need to solve (2) for a sequence of $\gamma > 0$. So the stability of the performance of the optimization algorithm with different $\gamma$ is very important. In our experiments, we have found that the performance of AMA is rather sensitive to the value of $\gamma$ in that the times taken to solve problems with different values of $\gamma$ can vary widely. However, SSNAL is much more stable. In Figure 2, we show the comparison between SSNAL and AMA on both the Two Half Moon and MNIST datasets with $\gamma \in [0.2 : 0.2 : 10]$.

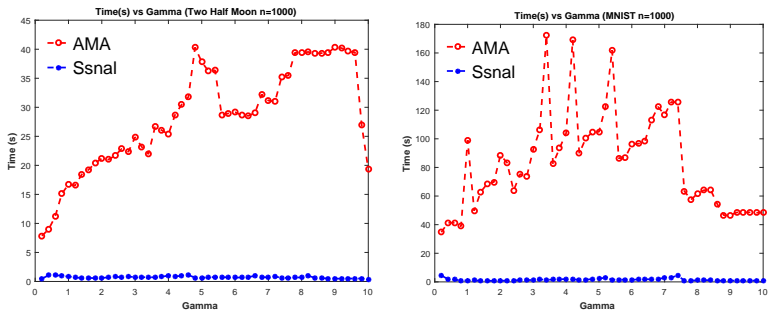

Figure 2: Time Comparison results between SSNAL and AMA on both Two Half Moon and MNIST datasets with $\gamma \in [0.2 : 0.2 : 10]$.

### 4.4 Scalability of Our Proposed Algorithm

In this section, we demonstrate the scalability of our algorithm SSNAL. Before we show the numerical results, we give some insights as to why our algorithm could be scalable. Recall that, the most computationally expensive step in our framework is in using the semismooth Newton-CG method to solve (10). However, if we look inside the algorithm, we can see that the key step is to use the CG method to solve (11) efficiently to get the Newton direction. According to our complexity analysis in Section 3.5, the computational cost for one step of CG update is $O(d|\widehat{\mathcal{E}}|)$. By the specific meaning of $\widehat{\mathcal{E}}$ explained in Section 3.5, $|\widehat{\mathcal{E}}|$ should grow only slowly with the growth of $n$. This low computational cost for the matrix-vector product in our Newton-CG method is the key point

behind why our algorithm can be scalable and efficient. The numerical results shown in this section also strongly support this argument.

In our experiments, we apply our algorithm on Half Moon data with $n$ ranging from 100 to 20000. Comparing to the numerical experiment results reported in [2] and [9] with $n \leq 500$ and $n \leq 600$, respectively, our results have convincingly demonstrated the scalability of our SSNAL.

In our experiments, we set $\phi = 0.5$, $k = 10$ (the number of nearest neighbors). Then we solve (2) with $\gamma \in [0.4 : 0.4 : 20]$. After generating the clustering path, we compute the average time for solving a single instance of (2) for each problem scale. Another factor related to the scalability is the number of neighbors $k$ used in $\mathcal{E}$ in (2). So, we also show the performance of SSNAL with different values of $k$. For each $k \in [5 : 5 : 50]$, we generate the clustering path for Two Half Moon dataset with $n = 2000$. Then we report the average time for solving a single instance of (2) for each $k$. We summarize our numerical results in Figure 3. We can observe that the computation time grows almost linearly with $n$ and $k$.

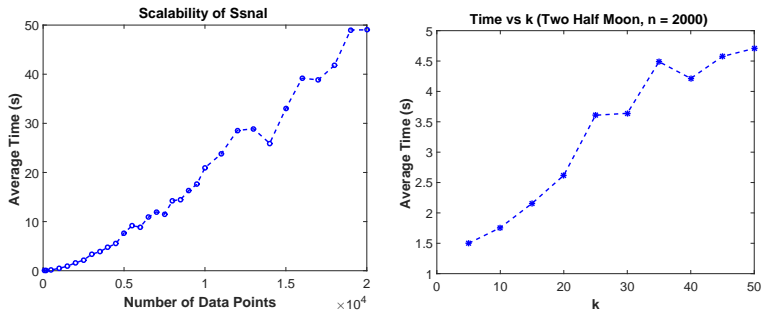

Figure 3: Numerical results to demonstrate the scalability of our proposed algorithm SSNAL with respect to $n$ and $k$.

## 5 Conclusion

In this paper, we proposed a highly efficient and scalable semismooth Newton based augmented Lagrangian method to solve the convex clustering model (2). To the best of our knowledge, this is the first optimization algorithm for convex clustering model which uses the second-order generalized Hessian information. Extensive numerical results shown in the paper have demonstrated the scalability and superior performance of our proposed algorithm SSNAL comparing to the state-of-the-art software CVXCLUSTR. The convergence results for our algorithm are also provided. In our future work, we plan to design a distributed and parallel version of SSNAL with the aim to handle huge scale datasets. From the modeling perspective, we will also work on generalizing our algorithm to handle kernel based convex clustering models.

## References

[1] L. CHEN, D. SUN, AND K. TOH, *An efficient inexact symmetric Gauss–Seidel based majorized*

*ADMM for high-dimensional convex composite conic programming*, Mathematical Programming, 161 (2017), pp. 237–270.

[2] E. Chi and K. Lange, *Splitting methods for convex clustering*, J. Computational and Graphical Statistics, 24 (2015), pp. 994–1013.

[3] F. Clarke, *Optimization and Nonsmooth Analysis*, John Wiley and Sons, New York, 1983.

[4] J.-B. Hiriart-Urruty, J.-J. Strodiot, and V. Nguyen, *Generalized Hessian matrix and second-order optimality conditions for problems with $C^{1,1}$ data*, Appl. Math. Optim., 11 (1984), pp. 43–56.

[5] T. D. Hocking, A. Joulin, F. Bach, and J.-P. Vert, *Clusterpath an algorithm for clustering using convex fusion penalties*, in 28th International Conference on Machine Learning, 2011.

[6] X. Li, D. Sun, and K. Toh, *A highly efficient semismooth Newton augmented Lagrangian method for solving lasso problems*, SIAM J. Optimization, ((in print)).

[7] F. Lindsten, H. Ohlsson, and L. Ljung, *Clustering using sum-of-norms regularization: With application to particle filter output computation*, in Statistical Signal Processing Workshop (SSP), IEEE, 2011, pp. 201–204.

[8] D. G. Mixon, S. Villar, and R. Ward, *Clustering subgaussian mixtures by semidefinite programming*, arXiv preprint arXiv:1602.06612, (2016).

[9] A. Panahi, D. Dubhashi, F. Johansson, and C. Bhattacharyya, *Clustering by sum of norms: Stochastic incremental algorithm, convergence and cluster recovery*, in 34th International Conference on Machine Learning, vol. 70, PMLR, 2017, pp. 2769–2777.

[10] M. Rezaei and P. Fränti, *Set-matching methods for external cluster validity*, IEEE Trans. on Knowledge and Data Engineering, 28 (2016), pp. 2173–2186.

[11] R. Rockafellar, *Augmented Lagrangians and applications of the proximal point algorithm in convex programming*, Mathematics of Operations Research, 1 (1976), pp. 97–116.

[12] ———, *Monotone operators and the proximal point algorithm*, SIAM J. Control and Optimization, 14 (1976), pp. 877–898.

[13] K. M. Tan and D. Witten, *Statistical properties of convex clustering*, Electronic J. Statistics, 9 (2015), p. 2324.

[14] C. Zhu, H. Xu, C. Leng, and S. Yan, *Convex optimization procedure for clustering: Theoretical revisit*, in Advances in Neural Information Processing Systems 27, 2014, pp. 1619–1627.